\documentclass[12pt]{article}
\input{xypic}
\usepackage{amsmath}
\usepackage{amssymb}
\usepackage{latexsym}
\usepackage{enumerate}

\renewcommand{\baselinestretch}{1.3}

\newtheorem{theorem}{Theorem}[section]
\newtheorem{proposition}[theorem]{Proposition}
\newtheorem{corollary}[theorem]{Corollary}
\newtheorem{lemma}[theorem]{Lemma}
\newtheorem{example}[theorem]{Example}
\newtheorem{remark}[theorem]{Remark}

\newenvironment{proof}{\noindent{\sc Proof.}}{\hfill\qed}

\newcommand{\Z}{{\mathbb Z}}

\newcommand{\C}{{\mathbb C}}

\newcommand{\Hom}{{\rm Hom}}
\newcommand{\Aut}{{\rm Aut}}
\newcommand{\Out}{{\rm Out}}

\newcommand{\Ker}{{\rm Ker}\,}

\newcommand{\qed}{\quad\lower0.05cm\hbox{$\Box$}}

\newcommand{\AG}{{\Aut(G)}}
\newcommand{\AH}{{\Aut(H)}}
\newcommand{\OG}{{\Out(G)}}
\newcommand{\OH}{{\Out(H)}}

\newcommand{\ls}[2]{{\,^{#1}\!#2}}                    

\newcommand{\downarrowright}[1]{\downarrow
\rlap{\raise0.1cm\hbox{$\scriptstyle{#1}$}}}

\newcommand{\downarrowleft}[1]{\rlap{\kern-0.2cm
\raise0.1cm\hbox{$\scriptstyle{#1}$}}\downarrow}

\newcommand{\uparrowright}[1]{\uparrow
\rlap{\lower0.1cm\hbox{$\scriptstyle{#1}$}}}

\newcommand{\uparrowleft}[1]{\rlap{\kern-0.2cm
\lower0.1cm\hbox{$\scriptstyle{#1}$}}\uparrow}

\newcommand{\ra}{\rightarrow}
\newcommand{\lra}{\longrightarrow}

\newcommand{\epi}{\mbox{$\to$\hspace{-0.35cm}$\to$}}

\newcommand{\mono}{\hookrightarrow}

\def\rmono{\rto|<\hole|<<\ahook}

\def\umono{\ar@{_{(}->}[u]}
\def\uumono{\ar@{_{(}->}[uu]}

\def\lmono{\ar@{_{(}->}[l]}
\def\llmono{\ar@{_{(}->}[ll]}

\begin{document}
\title{Non-simple localizations of finite simple groups}

\author{
{\sc Jos\'e L. Rodr\'{\i}guez, J\'er\^ome Scherer, and Antonio Viruel}
\thanks{The first and third authors were partially supported by
EC grant HPRN-CT-1999-00119 and CEC-JA grant FQM-213, the first and second authors
were partially supported by the Swiss National Science Foundation, the first author by
DGIMCYT grant BFM2001-2031, and the third author by MCYT-FEDER grant BFM2001-1825.}}

\date{}
\maketitle


\begin{abstract}
Often a localization functor (in the category of groups) sends a
finite simple group to another finite simple group. We study when
such a localization also induces a localization between the
automorphism groups and between the universal central extensions.
As a consequence we exhibit many examples of localizations of
finite simple groups which are not simple.
\end{abstract}


\section*{Introduction}

A group homomorphism $\varphi: H\to G$ is said to be a {\em
localization} if and only if $\varphi$ induces a bijection
\addtocounter{theorem}{1}
\begin{equation}
\label{bijection}
\varphi^*: \Hom(G,G)\cong \Hom(H, G)
\end{equation}
This is an ad hoc definition which comes from \cite[Lemma
2.1]{Cas00}. More details on localizations can be found there or
in the introduction of~\cite{RST}, where we exclusively study
localizations $H \mono G$, with both $H$ and $G$ simple groups.
Due to the tight links with homotopical localizations much effort
has been dedicated to analyze which algebraic properties are
preserved under localization. An exhaustive survey about this
problem is nicely exposed in \cite{Cas00} by Casacuberta. For
example, if $H$ is abelian and $\varphi: H \to G$ is a
localization, then $G$ is again abelian. Similarly, nilpotent
groups of class at most~2 are preserved (see~\cite[Theorem
3.3]{Lib99}), but the question remains open for arbitrary
nilpotent groups. Finiteness is not preserved, as shown by the
example $A_n\to SO(n-1)$ (this is the main result in
\cite{Lib98}). In the present paper we focus on simplicity of
finite groups and answer negatively a question posed both by
Libman in \cite{Lib99} and Casacuberta in \cite{Cas00} about
preservation of simplicity. In these papers it was also asked
whether perfectness is preserved. This is not the case either, as
we show with totally different methods in \cite{RSV}.

Our main result here is that if $H\mono G$ is a localization with
$H$ simple then $G$ need not be simple in general, see
Corollary~\ref{notsimple}. There is for example a localization map
from the Mathieu group $M_{11}$ to the double cover of the Mathieu
group $M_{12}$. This is achieved by a thorough analysis of the
effect of a localization on the Schur multiplier, which encodes
the information about the universal central extension. More
precisely we prove the following:

\medskip

\bigskip\noindent
{\bf Theorem 1.5} {\it
Let $i: H\mono G$ be an inclusion of two non-abelian finite simple groups and
$j: \tilde H \to \tilde G$ be the induced homomorphism on the universal
central extensions. Assume that $G$ does not contain any non-trivial central extension
of $H$ as a subgroup. Then $\,i: H\mono G$ is a localization if and only if
$\,j: \tilde H \to \tilde G$ is  a localization.
}

\medskip

We only consider non-abelian finite simple groups since the localization
of a cyclic group of prime order is either trivial or itself (\cite[Theorem~3.1]{Cas00}).
Naturally the second part of the paper deals with the effect of a localization
on the outer automorphism group, which roughly speaking is dual to the
Schur multiplier as it encodes the information about the ``super-group" of all
automorphisms.

\medskip

\bigskip\noindent
{\bf Theorem 2.4} {\it
Let $i:H\mono G$ be a localization between two non-abelian finite simple groups.
It extends then to a monomorphism $j: \AH \mono \AG$, which we assume induces
an isomorphism $\OH \cong \OG$. Then $j: \AH \mono \AG$ is a localization.}

\medskip

The converse does not hold: There exists a localization $\Aut(L_3(2)) \mono S_8$,
but the induced morphism $L_3(2) \mono A_8$ fails to be one.

\bigskip

{\it Acknowledgments:} We would like to thank Jon Berrick and Jacques Th\'evenaz
for helpful comments.

\section{Preservation of simplicity}
\label{section counter-example}

We first need to fix some notation. Let $Mult(G) = H_2(G; \Z)
\cong H^2(G; \C^{\times})$ denote the {\em Schur multiplier} of
the finite simple group $G$ and $Mult(G) \mono \tilde G \epi G$ be
its universal central extension. In particular the only
non-trivial endomorphisms of $\tilde G$ are automorphisms. This is
due to the fact that the only proper normal subgroups of $\tilde
G$ are contained in $Mult(G)$ and $Hom(G, \tilde G) = 0$ since the
universal central extension is not split. For more details, a good
reference is \cite[Section 6.9]{95f:18001}. Recall also that a
group $G$ is {\em perfect} if it is equal to its commutator
subgroup. Equivalently $G$ is perfect if $H_1(G; \Z) = 0$. If
moreover $H_2(G; \Z) = 0$ we say that $G$ is {\em superperfect}.
Hence for a perfect group $G$ we have that $\tilde G = G$ if and
only if $G$ is superperfect.

\medskip

Is simplicity preserved under localization?
We next show that the answer is affirmative if $H$ is maximal in $G$.
By $C_p$ we denote a cyclic group of order~$p$.

\begin{proposition}
\label{allissimple}
Let $G$ be a finite group and let $H$ be a maximal subgroup which is simple.
If the inclusion $H\mono G$ is a localization, then $G$ is simple.
\end{proposition}

\begin{proof}
First notice that $H$ cannot be normal in $G$. Indeed if $H$ is normal, the maximality of~$H$
implies that the quotient $G/H$ does not have any non-trivial proper subgroup.
Hence $G/H\cong C_p\,$ for some prime~$p$. But then $G$ has a subgroup of order~$p$
and there is an endomorphism of $G$ factoring through $C_p\,$, whose restriction to
$H$ is trivial. This contradicts the assumption that the inclusion $H\mono
G$ is a localization.

Let $N$ be a normal subgroup of $G$. As $H$ is simple, $N\cap H$ is either equal to
$\{1\}$ or~$H$. If $N\cap H=H$, as $H$ is maximal, then either $N=G$ or $N=H$, and we
just showed that the latter case is impossible.
If $N\cap H=\{1\}$, then either $N=\{1\}$ or $NH=G$ as $H$ is maximal. The
second case cannot occur because it would imply that $G=N\rtimes H$, but
$H\mono N\rtimes H$ cannot be a localization since both the identity
of~$G$ and the projection onto~$H$ extend the inclusion ${H\mono G}$. Therefore
there are no normal proper non-trivial subgroups in $G$.
\end{proof}

\bigskip

We indicate next (in Corollary~\ref{notsimple}) a generic situation where the localization
of a simple group can be non-simple (it will actually be the universal cover of a simple
group). To achieve this we study when a localization of finite simple groups induces
a localization of the universal covers.

\begin{proposition}
\label{help}
Let $H$ and $G$ be non-abelian finite simple groups. Assume that any
homomorphism between the universal central extensions $\tilde H \to \tilde G$
sends $Mult(H)$ to $Mult(G)$. Then $p: \tilde G \epi G$ and $q: \tilde H \epi H$ induce
an isomorphism $F: \Hom(\tilde H, \tilde G) \stackrel{\simeq}{\lra} \Hom(H, G)$.
\end{proposition}

\begin{proof}
First notice that $p$ and $q$ induce indeed a map
$F: \Hom(\tilde H, \tilde G) \to \Hom(H, G)$ by our assumption that any
homomorphism $\tilde H \to \tilde G$ sends the center to the center.
We show now that $F$ is surjective.
Let $\alpha: H \to G$. Using the $k$-invariants $k_H:K(H, 1) \ra K(Mult(H), 2)$
and $k_G:K(G, 1) \ra K(Mult(G), 2)$ classifying the
universal central extensions, construct the commutative diagram
\[
\diagram
K(H, 1) \rto^\alpha \dto_{k_H} & K(G,1) \dto^{k_G} \cr
K(Mult(H), 2) \rto_{H_2(\alpha)_*} & K(Mult(G), 2)
\enddiagram
\]
Taking vertical fibres gives precisely a map $K(\tilde H, 1) \ra K(\tilde G, 1)$
induced by some morphism $\beta: \tilde H \ra \tilde G$ with $F(\beta) = \alpha$.
Let us show now that $F$ is also injective by indicating an equivalent construction.
Given a morphism $\alpha: H \ra G$, construct
the pull-back $P_\alpha$ of $p$ along $\alpha$. Then $P_\alpha \epi H$ is
a central extension, so that there exists a unique compatible morphism
$\tilde H \ra P_\alpha$. The composite $\tilde H \ra P_\alpha \ra \tilde G$ is hence
the unique morphism whose image under $F$ is $\alpha$.
\end{proof}

\begin{corollary}
\label{aut}
Let $G$ be a non-abelian finite simple group and denote by
$p: \tilde G \epi G$ its universal
central extension. Then we have an isomorphism
$F: \Aut(\tilde G) \stackrel{\simeq}{\lra} \Aut(G)$.
\end{corollary}

\begin{proof}
We have to check that any homomorphism $\tilde G \ra \tilde G$ sends
the center to the center. As the only morphism which is not an automorphism
is the trivial one, this is a clear consequence of the fact that the
image of the center is contained in the center of the image.
The proposition tells us that we have an isomorphism
$F: \Hom(\tilde G, \tilde G) \stackrel{\simeq}{\lra} \Hom(G, G)$, therefore also
one $F: \Aut(\tilde G) \stackrel{\simeq}{\lra} \Aut(G)$.
\end{proof}

\medskip

One should be warned that this result does not imply that an automorphism of the
universal central extension always induce the identity on the center (of course
all inner automorphisms do so). For example let $G = L_3(7) = A_2(7)$, so
$\tilde L_3(7) = SL_3(7)$ and $Mult(L_3(7)) = Z(SL_3(7)) \cong \Z/3$ is
generated by the diagonal matrix $D$ whose coefficients are $2$'s. There is
an outer ``graph automorphism" of order $2$ given by the transpose of the
inverse. It sends a matrix $A$ to $\ls{t}{A^{-1}}$, so the image of $D$ is $D^{-1}$.

\begin{proposition}
Let $G$ be a finite simple group.
Then, the universal cover $\tilde G\epi G$ is a localization.
\end{proposition}

\begin{proof}
We have to show that $\tilde G\epi G$ induces a bijection $\Hom(G,
G) \cong \Hom(\tilde G,G)$ or equivalently, $\Aut(G) \cong
\Hom(\tilde G,G) \setminus \{0\}$. This follows easily since the
only non-trivial proper normal subgroups of $\tilde G$ are
contained in its center $Mult(G)$. Thus any non-trivial
homomorphism $\tilde G\to G$ can be decomposed as the canonical
projection $\tilde G\to G$ followed by an automorphism of~$G$.
\end{proof}

\medskip

\begin{theorem}
\label{covers}
Let $i: H\mono G$ be an inclusion of two non-abelian finite simple groups and
$j: \tilde H \to \tilde G$ be the induced homomorphism on the universal
central extensions. Assume that $G$ does not contain any non-trivial central extension
of $H$ as a subgroup. Then $\,i: H\mono G$ is a localization if and only if
$\,j: \tilde H \to \tilde G$ is  a localization.
\end{theorem}

\begin{proof}
The map $i: H\mono G$ is a localization if and only if it induces an isomorphism
$\Hom(H, G) \cong \Hom(G, G)$. Let us analyze the behavior of morphisms
$\varphi: \tilde H \ra \tilde G$. By composing with $q: \tilde G \to G$ we get
a morphism $\tilde H \ra G$. As $G$ does not contain any subgroup isomorphic
to a central extension of $H$, we see that $\varphi(Mult(H)) \subset Mult(G)$.
We deduce now by Proposition~\ref{help} that the
universal central extensions induce isomorphisms
$\Hom(\tilde H, \tilde G) \cong \Hom(H, G)$ as well as
$\Hom(\tilde G, \tilde G) \cong \Hom(G, G)$. Both isomorphisms
are compatible, so $i: H\mono G$ is a localization if and only if
$\,j: \tilde H \to \tilde G$ induces an isomorphism
$\Hom(\tilde H, \tilde G) \cong \Hom(\tilde G, \tilde G)$.
\end{proof}

\medskip

\begin{remark}
\label{centers}
{\rm We do not know how to remove the assumption on the centers in
Proposition~\ref{help}. There exist indeed morphisms between
covers of finite simple groups which do not send the center to the
center. One example is given in \cite[p.34]{atlas} by the
inclusion $\tilde A_5 \mono U_3(5)$. A larger class of examples is
obtained as follows: Let $H$ be a finite simple group of order $k$
and $\tilde H$ its universal central extension of order
$n=|Mult(H)| \cdot k$. The regular representations $H \mono S_k$
and $\tilde H \mono S_n$ actually lie in $A_k$ and $A_n$ because
the groups are perfect. Therefore $A_n$ contains both $H$ and
$\tilde H$ as subgroups. However we do not know of a single
example of a localization $H \mono G$ which does not satisfy this
assumption and it is rather easy to check in practice. }
\end{remark}

\noindent {\bf Question}: Let $i: H \mono G$ be a localization. Is
it possible that some subgroup of $G$ be isomorphic to a
non-trivial central extension of $H$? If the answer is no, we
would get a more general version of Theorem~\ref{covers}. This
would form a perfectly dual result to Theorem~\ref{autlocalize} if
the extra assumption that $i$ induces an isomorphism $H_2(i; \Z):
Mult(H) \ra Mult(G)$ has to be used.

\medskip

Beware that in general the induced morphism on the universal
central extensions given by the above theorem is not an inclusion.
For example $L_2(11) \mono U_5(2)$ is a localization by the main
theorem in~\cite{RST}. However $U_5(2)$ is superperfect and the
universal central extension $SL_2(11)$ of $L_2(11)$ is not a
subgroup of $U_5(2)$. Nevertheless there is a localization
$SL_2(11) \to U_5(2)$. The dual situation when $H$ is superperfect
leads to our counterexamples.

\begin{corollary}
\label{notsimple}
Let $i: H\mono G$ be an inclusion of two non-abelian finite simple groups and
assume that $H$ is superperfect. Let also
$j: H = \tilde H \mono \tilde G$ denote the induced homomorphism on the universal
central extensions. Then $\,i: H\mono G$ is a localization if and only if
$\,j: H \mono \tilde G$ is  a localization.
\end{corollary}

\begin{proof}
There are no non-trivial central extensions of $H$ so Theorem~\ref{covers} applies.
\end{proof}

\begin{example}
\label{viru}
{\rm
The inclusion $M_{11}\mono \tilde M_{12}$ of the Mathieu group~$M_{11}$
into the double cover of the Mathieu group~$M_{12}$ is a localization.
This follows from the above proposition.
Note that $M_{11}$ is not maximal in $\tilde M_{12}$
(the maximal subgroup is $M_{11}\times C_2$), so this does not contradict
Proposition~\ref{allissimple}. The following inclusions are localizations:
$Co_2\mono Co_1$ and $Co_3\mono Co_1$ by \cite[Section 4]{RST}.
As the smaller group is superperfect we get localizations
$Co_2\mono \tilde Co_1$ and $Co_3\mono \tilde Co_1$.

We get many other examples of this type
using \cite[Corollary 2.2]{RST}. All sporadic groups appearing in this corollary
which have trivial Schur multiplier (that is $M_{11}$, $M_{23}$, $M_{24}$, $J_1$, $J_4$,
$Co_2$, $Co_3$, $He$, $Fi_{23}$, $HN$, and $Ly$) admit the double cover of an alternating
group as localization (as $Mult(A_n)$ is cyclic of order 2 for $n > 7$).
}
\end{example}

\begin{remark}
\label{monster}
{\rm
The inclusion $Fi_{23} \mono B$ of the Fischer group into the baby monster is a localization
by \cite[Section 3 (vi)]{RST}. This yields a localization $Fi_{23} \mono \tilde B$.
As the double cover $\tilde B$ is a maximal subgroup of the Monster~$M$,
it would be nice to know if $\tilde B \mono M$ is a localization. This would connect
the Monster to the rigid component of the alternating groups (in \cite{RST} we
were able to connect all other sporadic groups to an alternating group by a zigzag
of localizations).
}
\end{remark}

\section{Localizations between automorphism groups}
\label{section auto}

The purpose of this section is to show that a localization $H \mono G$
can often be extended to a localization $\Aut(H) \mono \Aut(G)$, similarly
to the dual phenomenon observed in Theorem~\ref{covers}.
This generalizes the observation made by Libman
(cf. \cite[Example 3.4]{Lib99}) that the localization $A_n \mono A_{n+1}$
extends to a localization $S_n \mono S_{n+1}$ if $n\geq 7$.
This result could be the starting point for determining
the rigid component (as defined in \cite{RST}) of the symmetric groups, but
we will not go further in this direction.

\begin{lemma}
\label{normalsubgroup}
Let $G$ be a non-abelian finite simple group. Then any proper normal subgroup of
$\AG$ contains $G$. In particular any endomorphism of $\AG$ is either an isomorphism,
or contains $G$ in its kernel.
\end{lemma}

\begin{proof}
Let $N$ be a normal subgroup of $\AG$ and assume that it does not contain~$G$. Since
$N \cap G$ is a normal subgroup of $G$, it must be the trivial subgroup. Hence the
composite $N \mono \AG \epi \OG$ is injective. The orbit under conjugation by $G$ of an
automorphism in $\AG$ is reduced to a point if and only if the automorphism is the
identity. Thus $N$ has to be trivial.
\end{proof}

\begin{lemma}
\label{Hinaut}
Let $G$ be a non-abelian finite simple group. Then any non-abelian simple subgroup of $\AG$
is contained in $G$.
\end{lemma}

\begin{proof}
Let $H$ be a non-abelian simple subgroup of $\AG$. The kernel $G$ of the projection
$\AG \epi \Out(G)$ contains $H$ because $\Out(G)$ is
solvable (this is the Schreier conjecture, whose proof depends on the
classification of finite simple groups, see \cite[Theorem~7.1.1]{98j:20011}).
\end{proof}

\begin{lemma}
\label{unique}
Let $i:H\mono G$ be a localization between two non-abelian finite simple groups.
Then it extends in a unique way to a monomorphism $j: \AH \mono \AG$.
\end{lemma}

\begin{proof}
Since $i$ is a localization, it extends to a monomorphism $j: \AH \mono \AG$
by Theorem~1.4 in \cite{RST}. The uniqueness is given by \cite[Remark~1.3]{RST}.
Indeed, given a commutative square
\[
\diagram
H \rmono^i \dto & G \dto \cr
\AH \rmono_j & \AG
\enddiagram
\]
Lemma~1.2 in \cite{RST} implies that, for any $\alpha \in \AH$, $j(\alpha)$ is an automorphism of
$G$ extending $\alpha$. But there exists a unique automorphism $\beta: G \ra G$ such that
$\beta \circ i = i \circ \alpha$ because $i$ is a localization.
\end{proof}

\begin{theorem}
\label{autlocalize}
Let $i:H\mono G$ be a localization between two non-abelian finite simple groups.
It extends then to a monomorphism $j: \AH \mono \AG$, which we assume induces
an isomorphism $\OH \cong \OG$. Then $j: \AH \mono \AG$ is a localization.
\end{theorem}

\begin{proof}
As $i$ is a localization it extends to a unique inclusion $j: \AH
\mono \AG$ by the above lemma. Let $\varphi: \AH \ra \AG$ be any
homomorphism. We have to show that there is a unique $\alpha: \AG
\ra \AG$ such that $\alpha \circ i = \varphi$. If $\alpha: \AG \to
\AG$ is not an isomorphism, it factorizes through some quotient
$Q$ of $\OG$ by Lemma~\ref{normalsubgroup}. The assumption that
$j$ induces an isomorphism on the outer automorphism groups
implies then that the restriction of $\alpha$ to $\AH$ is trivial
if and only if $\alpha$ is trivial. Therefore if $\varphi$ is
trivial, we conclude that the unique such $\alpha$ is the trivial
homomorphism.

Let us assume that $\varphi$ is not trivial. If it is an
injection, the image of the composite $\psi: H \mono \AH
\stackrel{\varphi}{\lra} \AG$ actually lies in $G$ by
Lemma~\ref{Hinaut} and because $H \mono G$ is a localization,
there is a unique automorphism $\alpha$ of $G$ making the
appropriate diagram commute. Conjugation by $\alpha$ on $\AG$ is
the unique extension we need. Indeed in the following diagram all
squares are commutative and so is the top triangle:
\[
\diagram
H \xto[rr]^i \xto[dr]^\psi \xto[dd] && G \xto[dd] \xto[dl]_\alpha \cr
& G \xto[dd] & \cr
\AH \xto[rr]|(0.50)\hole^(0.60){i} \xto[dr]^\varphi && \AG
\xto[dl]_{c_\alpha} \cr
& \AG &
\enddiagram
\]
The map $\psi$ is also a localization so that
the bottom triangle commutes as well by Lemma~\ref{unique}. If $\varphi$ is not injective,
then $H \subset \Ker \varphi$ by Lemma \ref{normalsubgroup}.
In that case the image of $\varphi$ in $\AG$ is some quotient of $\OH$.
As $j$ induces an isomorphism $\OH \cong \OG$, $\varphi$ clearly extends
to a unique endomorphism of $\AG$.
\end{proof}

\medskip

The assumption that $j$ induces the isomorphism between the outer
automorphism groups cannot be dropped. In fact even when $\OH$ and
$\OG$ are cyclic of order~$2$, a localization $H \mono G$ does not
always induce one $\AH \mono \AG$. For example $i: L_3(3) \mono
G_2(3)$ is a localization (see \cite[Proposition~4.2]{RST}), but
$\Aut(L_3(3))$ is actually contained in $G_2(3)$. Thus $j:
\Aut(L_3(3)) \mono \Aut(G_2(3))$ cannot be a localization for the
good and simple reason that the composite $\Aut(L_3(3)) \mono
\Aut(G_2(3)) \epi \Out(G_2(3)) \cong C_2 \mono \Aut(G_2(3))$ is
trivial. The same phenomenon occurs again for $i: He \mono
Fi_{24}'$. Still, many examples can be directly derived from this
theorem, as it is often easy to check that $j$ must induce an
isomorphism $\OH \cong \OG$.

\begin{corollary}
\label{corautlocalize}
Let $i:H\mono G$ be a localization between two non-abelian finite simple groups.
Assume that $H$ is a maximal subgroup of $G$ and that both $\OH$ and $\OG$
are cyclic groups of order $p$ for some prime $p$.
Then $j: \AH \mono \AG$ is a localization.
\end{corollary}

\begin{proof}
We only have to show that $j$ itself induces the isomorphism $\OH \cong \OG$.
Because these outer automorphism groups are cyclic of prime order, the induced
morphism must be either trivial or an isomorphism. Being trivial means that
any automorphism of $H$ is sent by $j$ to an inner automorphism of $G$, which
can happen only if $\Aut(H)$ is a subgroup of $G$.
\end{proof}

\medskip

Directly from the corollary we deduce that $S_n  \mono S_{n+1}$ and
$SL_2(p) \mono S_{p+1}$ are localizations (by \cite[Proposition~2.3(i)]{RST}
$L_2(p) \mono A_{p+1}$ is a localization). Suzuki's
chain of groups $L_2(7) \mono G_2(2)' \mono J_2 \mono G_2(4) \mono Suz$
(see \cite[p.108-9]{gore})
also extends to localizations of their automorphism
groups
$$
\Aut(L_2(7)) \mono \Aut(G_2(2)') \mono \Aut(J_2) \mono \Aut(G_2(4)) \mono \Aut(Suz).
$$

\medskip

\begin{remark}
\label{falseconverse}
{\rm
The converse of the above theorem is false.
There exists for example an inclusion $\Aut(L_3(2)) \mono S_8$ which
is actually a localization (Condition (\ref{bijection}) can be checked quickly
with the help of MAGMA). However the induced morphism $L_3(2) \mono
A_8$ fails to be a localization: There are two conjugacy classes of subgroups
of $A_8$ isomorphic to $L_3(2)$, which are not conjugate in $S_8$.
}
\end{remark}


\vskip 0.5 cm

\setlength{\baselineskip}{0.6cm}

\bigbreak

Jos\'e L. Rodr\'{\i}guez

\noindent
Departamento de Geometr\'\i a, Topolog\'\i a y Qu\'\i mica Org\'anica,
Universidad de Almer\'\i a, E--04120 Almer\'\i a, Spain,
e-mail: {\tt jlrodri@ual.es}

\bigskip

J\'er\^ome Scherer

\noindent Departament de Matem\`atiques, Universitat Aut\`onoma de
Barcelona, E--08193 Bellaterra, Spain, e-mail: {\tt
jscherer@mat.uab.es}

\bigskip

Antonio Viruel

\noindent
Departamento de Geometr\'\i a y Topolog\'\i a,
Universidad de M\'alaga, E--04120 M\'alaga, Spain,
e-mail: {\tt viruel@agt.cie.uma.es}

\end{document}